\documentclass{amsart}

\usepackage{amsmath,amsthm,amssymb}

\usepackage{pb-diagram,pstricks,pst-tree}

\newcommand{\id}{\operatorname{Id}}

\newcommand{\ad}{\operatorname{ad}}

\newcommand{\Ad}{\operatorname{Ad}}

\renewcommand{\ker}{\operatorname{Ker}}

\newcommand{\im}{\operatorname{Im}}

\newcommand{\coder}{\operatorname{Coder}}

\newcommand{\comap}{\operatorname{Comap}}

\newcommand{\Def}{\operatorname{\mathit{Def}}}

\newcommand{\g}{\mathfrak{g}}

\newcommand{\h}{\mathfrak{h}}
\newcommand{\M}{\mathcal{M}}

\renewcommand{\O}{\mathcal{O}}

\newcommand{\<}{\langle}
\renewcommand{\>}{\rangle}

\theoremstyle{plain}
\newtheorem{theorem}{Theorem}[section]
\newtheorem{proposition}[theorem]{Proposition}
\newtheorem{corollary}[theorem]{Corollary}
\theoremstyle{definition}
\newtheorem{definition}[theorem]{Definition}
\newtheorem{remark}[theorem]{Remark}
\newtheorem{example}[theorem]{Example}

\begin{document}
\title{Deformations of associative algebras with inner products}
\author{John Terilla}
\email{jterilla@math.sunysb.edu}
\address{Department of Mathematics\\ 
                State University of New York\\
                Stony Brook, NY\\
                11794-3651.}
  
\author{Thomas Tradler}
\email{tradler@stanford.edu}
\address{Department of Mathematics\\
        Bldg. 380
        Stanford University\\
        Stanford, CA 94305-2125.}
        
\thanks{The authors would like to thank Jim Stasheff, Martin Markl, and Dennis Sullivan for many helpful discussions.}


\begin{abstract}
We develop the deformation theory of $A_\infty$ algebras together
with $\infty$-inner products and identify a differential graded
Lie algebra that controls the theory.  This generalizes the deformation
theories of associative algebras, $A_\infty$ algebras, associative
algebras with inner products, and $A_\infty$ algebras with inner
products.
\end{abstract}

\maketitle

\tableofcontents

\section{Introduction}

In this paper, we describe a
differential graded Lie algebra that controls deformations of the
homotopy versions of associative algebras together with the inner
products introduced in \cite{T1}.  Let us review the basic idea of a deformation theory governed
by a differential graded Lie algebra \cite{SS2,GM,K1,M}.

Fix a ground field $k$ of characteristic $\neq2$.  For any differential graded Lie (dg Lie)
algebra $(\g=\oplus_i \g^i,d,[\,,\,])$ over $k$,  one can consider
deforming the differential $d$ in the direction of an inner derivation.
Informally, such a deformation is given by an (equivalence classes of) $\alpha$ making
$$d_\alpha : = d+\mathrm{ad}(\alpha)$$
into a differential.  The map $d_\alpha$ is always a derivation
and the condition that $d_\alpha^2=0$ translates into the
Mauer-Cartan equation:
$$d\alpha+\frac{1}{2}[\alpha,\alpha]=0.$$
The deformed differential $d_\alpha$ may involve parameters from
the maximal ideal $m$ of a $\mathbb{Z}$ graded Artin local ring:
$\alpha \in (\g \otimes_k m)^1$. If $m$ is the maximal ideal of a
local Artin ring $R$ and $\alpha\in (\g \otimes_k m)^1$ is a
solution to the Mauer-Cartan equation, then one may call
$d_\alpha$ a \emph{deformation of $d$ over $R$}.
Two deformations of $d$ should be considered
equivalent if they are related via a natural automorphism of $\g
\otimes m$.  A ring map $R \to S$ will transport a deformation of $d$ over $R$ to a deformation of $d$ over $S$.

More formally, one has a functor $\Def_\g$ from the category of
$\mathbb{Z}$ graded Artin local rings with residue field $k$ to the category of sets, assigning to such a ring $R$ with maximal ideal $m$
the set
$$\Def_\g(R)=\{\alpha \in (\g\otimes_k m)^1:
d\alpha+\frac{1}{2}[\alpha,\alpha]=0\}/\sim.$$ Here, $\sim$ is the
equivalence relation determined by the action of the gauge group,
which we now recall.  Since $R$ is an Artin ring, $m$ is a
nilpotent algebra, and $(\g\otimes_k m)^0 \subseteq
\g \otimes_k m$ is a nilpotent Lie algebra.  Therefore, there
exists a group $G=\{\exp{\beta}: \beta \in (\g \otimes_k m)^0\}$, called the gauge group,
with multiplication defined by the Baker-Campbell-Hausdorff
formula. The action of $e^\beta\in G$ on an element $\alpha \in
(g\otimes_k m)^1$ is determined by the infinitesimal action:
$$\alpha \mapsto \beta \cdot \alpha =
[\beta,\alpha]-d\beta, \quad \alpha \in (\g\otimes m)^1, \, \beta
\in (\g\otimes m)^0.$$
This action satisfies $$e^{\ad \beta} d_\alpha e^{-\ad
\beta}=d_{e^{\beta} \cdot \alpha},$$ and preserves the set of
solutions to the Maurer-Cartan equation. 
In this paper, we work with $A_\infty$ algebras equipped with an
$\infty$ inner product.  One has the notion of a deformation of an
$A_\infty$ algebra with an $\infty$ inner product over a
ring $R$, and there is a natural equivalence on the set of deformations.  
A ring map $R \to S$ transports deformations over $R$ to deformations over $S$.
The association
$$R \mapsto
\left\{
\begin{array}{c}
\text{deformations of the }A_\infty \text{ algebra with }\\ \text{the }\infty\text{  inner product over } R
\end{array}\right\} \biggl/
\left\{\begin{array}{c}
\text{equivalent}\\ \text{deformations}\end{array}
\right\}
$$
defines a covarient deformation functor.

In this paper,
we construct a dgLa $(\h=\oplus_i \h^i,d,[\,,\,])$
associated to an $A_\infty$ algebra with an $\infty$ inner product,
and prove that the functor described above is isomorphic to $\Def_\h$.
This is the precise mathematical content of the statement
\emph{the dgLa $(\h,[\,,\,],d)$ controls the deformations of the
$A_\infty$ algebra with an $\infty$ inner product.}

\section{Definitions of $A_\infty$ algebras and $\infty$ inner products}

We now review the concept of an $\infty$ inner product on an
$A_\infty$ algebra \cite{T1}, \cite{TZ}. The concepts of
$A_\infty$ algebras, $A_\infty$ bimodules, $A_\infty$ bimodule
maps, and $A_\infty$ inner products are generalizations of the
usual concepts of associative algebras, bimodules, bimodule maps,
and invariant inner products.

\subsection{$A_\infty$ algebras}
Let $V=\bigoplus_{j\in \mathbb{Z}} V^{j}$ be a graded module over
a ring $S$.  Recall that the suspension $V[1]$ of $V$ is defined
to be $V[1]= \bigoplus_{j\in \mathbb{Z}} (V[1])^{j}$ with
$(V[1])^{j}:= V^{j-1}$.  For a graded $S$-module $A$, $BA$ denotes
the tensor coalgebra $BA:=T(A[1])$, sometimes called the \emph{bar
complex} of $A$. An \emph{$A_\infty$ algebra over $S$} is defined
to be a pair $(A,D)$ where $A$ is a graded $S$ module and $D\in
\coder(BA)$ of degree $-1$ with $D^2 =0$.

Suppose that $(A,D)$ and $(A',D')$ are $A_\infty$ algebras over $S$.  Then an $A_\infty$ map from $(A',D')$ to $(A,D)$ is a map $\lambda:BA' \to BA$ satisfying $\lambda\circ D'=D\circ \lambda.$

\subsection{$A_\infty$ bimodules.}
Let $(A,D)$ be an $A_\infty$ algebra over $S$, and let $M$ be a
graded $S$ module. Let $B^{M}A$ denote the tensor bicomodule
$B^{M}A:=T^{M[1]}(A[1])=\bigoplus_{k,l\geq 0}A[1]^{\otimes k}
\otimes M[1] \otimes A[1]^{\otimes l}$ of $M[1]$ over $BA$. An
\emph{A$_\infty$ bimodule structure on $M$ over $A$} is defined to
be a coderivation $D^{M}\in \coder_{D}(B^{M}A,B^{M}A)$ over $D$ of
degree $-1$ with $(D^{M})^{2}=0$.


Let $(M,D^{M})$ and $(N,D^{N})$ be $A_\infty$ bimodules over $A$.
Let $\comap(B^{M}A,B^{N}A)$ denote the maps $F:B^{M}A\to
B^{N}A$ satisfying
$$
\begin{diagram}
  \node{B^M A}\arrow{e,t}{\Delta^{M}} \arrow{s,l}{F} \node{(BA\otimes
    B^M A)\oplus (B^M A\otimes BA)}\arrow{s,r}{(\id\otimes F)\oplus
    (F\otimes \id)}\\
  \node{B^N A}\arrow{e,b}{\Delta^{N}} \node{(BA\otimes B^M A)\oplus
    (B^M A\otimes BA)}
\end{diagram}
$$
The space $\comap(B^{M}A,B^{N}A)$ carries a differential defined by
$$\delta^{M,N}(F):=D^{N}\circ F-(-1)^{|F|}F\circ D^{M}.$$
In this case, an \emph{A$_\infty$ bimodule map from $M$ to $N$} is
defined to be an element $F\in \comap(B^{M}A,B^{N}A)$ of degree
$0$ with $\delta^{M,N}(F)=0$, i.e.
$$ D^{N}\circ F=F\circ D^{M}. $$

\subsection{$\infty$ inner products}
For any $f\in \coder(BA)$, there are induced coderivations
$f^{A}\in \coder_{f}(B^{A}A,B^{A}A)$ and $f^{A^{*}}\in \coder_{f}
(B^{A^{*}}A,B^{A^{*}}A)$, where $A^*=\hom_S(A,S)$ denotes the dual
of $A$.  One also has an induced map
$$\delta_f:\comap(B^{M}A,B^{N}A) \to \comap (B^{M}A,B^{N}A)$$
given by $\delta_f (F)=f^{A^*}\circ F- (-1)^{|f||F|}\cdot F\circ
f^A$.
Note that, in particular, if $(A,D)$ is an $A_\infty$ algebra, then $A$ and $A^*$ 
have $A_\infty$ bimodule structures given by $D^A$ and $D^{A^*}$.

\begin{definition}
Let $(A,D)$ be an $A_\infty$ algebra over $S$.  We define an
\emph{$\infty$ inner product on $A$ over $S$} to be an $A_\infty$
bimodule map $I$ from $A$ to $A^{*}$. Equivalently, an $\infty$ inner
product is an element $I\in \comap(B^{A}A,B^{A^{*}}A)$ satisfying
$$\delta_D(I)=D^{A^{*}}\circ I-I\circ D^{A}=0.$$
Every inner product $\<\,,\, \>:A\otimes A \to S$ defines an
element $I\in \comap(B^A A, B^{A^*} A)$. In this case, the
condition $D^{A^{*}}\circ I-I\circ D^{A}=0$ is equivalent to
$\<D(a_1,\ldots,a_n),a_{n+1}\>=\pm \<a_1,D(a_2,\ldots,
a_{n+1})\>$.  See the appendix for additional illustrations.
\end{definition}

\subsection{Induced maps}
Recall, that if $\lambda:A'\to A$ is an algebra map
between two associative algebras, then every module over
$A$ is also a module over $A'$, and similarly for module maps.
Also, $\lambda:A'\to A$ and $\lambda^*:A^*\to (A')^*$ will be module maps over $A'$. Here we give the corresponding homotopy generalizations.

Suppose that $\lambda$ is an $A_\infty$ map from $(A',D')$ to $(A,D)$.
First, every A$_\infty$ bimodule
$(M,D^M)$ over $A$ is also an $A_\infty$ bimodule over $A'$, whose
structure map is determined by the lowest components (which are maps $B^MA\to M$)
\begin{multline*}
(D^M)^{\lambda}(a'_1,\ldots ,a'_k,m,a'_{k+1},\ldots, a'_{k+l})\\
=\sum \pm pr_M\circ D^M (\lambda(a'_1,\ldots),\ldots,
\lambda(\ldots,a'_k),m, \lambda(a'_{k+1},\ldots),
\ldots,\lambda(\ldots,a'_{k+l})).
\end{multline*}
Here, $pr_M$ denotes the projection onto $M$.  
The signs are given by the usual sign rule, namely introducing a
sign $(-1)^{|\alpha| \cdot |\beta|}$, whenever $\alpha$ jumps over
$\beta$.  The relevant degrees are the degrees given
in $B^M A$.  

Also, any $A_\infty$ bimodule map $F:B^M A\to B^N A$ over $A$
induces an $A_\infty$ bimodule map $F^{\lambda}:B^M A'\to B^N A'$
over $A'$ given by
\begin{multline*}
F^{\lambda}(a'_1,\ldots,a'_k,m,a'_{k+1},\ldots, a'_{k+l})\\
=\sum \pm pr_N\circ F(\lambda(a'_1,\ldots),\ldots,
\lambda(\ldots,a'_k),m, \lambda(a'_{k+1},\ldots), \ldots,
\lambda(\ldots,a'_{k+l})).
\end{multline*}
Furthermore $\lambda$ induces the two A$_\infty$ bimodule maps
over $A'$
$$\overline{\lambda}:B^{A'} A' \to B^A A' \text{ and }
 \widetilde{\lambda}:B^{A^*} A' \to B^{(A')^*} A'$$
defined by the components
\begin{gather*}
\overline{\lambda}(a'_1,\ldots,a'_{k+l+1})=pr_A\circ
\lambda(a'_1,\ldots,a'_{k+l+1}) \intertext{ and }
 (\widetilde{\lambda}(a'_1,\ldots,a^*,\ldots,a'_{k+l}))(a')=\pm
  a^*(pr_A\circ \lambda(a'_{k+1},\ldots,a'_{k+l},a',a'_1,\ldots,a'_k)).
\end{gather*}

For our own comfort, we write $TA$ and $T^M A$ instead of $BA$ and
$B^M A$ for the remainder of this paper.

\section{Deformations of $A_\infty$ algebras and $\infty$ inner products}

Before we define the specific differential graded Lie algebra
$(\h,d,[\,,\,])$ that controls the deformations of $A_\infty$
structures and $\infty$ inner products, we discuss a simple
example and make a remark.

\begin{example}
Any graded associative algebra $\g$ becomes a Lie algebra by
defining the bracket to be the usual commutator. An element
$\alpha \in \g^1$ satisfying $\alpha^2=0$ is sometimes called a
\emph{polarization}.  With a polarization $\alpha\in \g^1$,
$\g$ becomes a differential graded Lie algebra
 by setting the differential to be $\delta=\ad(\alpha)$.
With $\delta$ so defined, the Mauer-Cartan equation becomes
$$ 0=\delta(\gamma)+\frac{1}{2}[\gamma,\gamma]=
   \frac{1}{2}[\alpha+\gamma, \alpha+\gamma].$$
In other words, $\gamma\in \g^1$ satisfies the Mauer-Cartan
equation if and only if $\alpha+\gamma$ is another polarization.

Now let $S$ be a graded ring and consider $\g$ defined by
$$ \g=\left\{%
\begin{pmatrix}
  a & 0 \\
  b & a \\
\end{pmatrix}%
\biggr | \, a,b \in S\right\}$$ with the bracket defined as the
usual graded commutator of matrix multiplication:
\begin{equation*}
 \left[
\begin{pmatrix}
  a & 0 \\
  b & a \\
\end{pmatrix},
\begin{pmatrix}
  c & 0 \\
  d & c \\
\end{pmatrix}
\right]=
\begin{pmatrix}
  [a,c] & 0 \\
  [b,c]+[a,d] & [a,c] \\
\end{pmatrix}.
\end{equation*}
Then,
$\begin{pmatrix}
  D & 0 \\
  I & D \\
\end{pmatrix} \in \g^1$
is a polarization if and only if
$$0=[D,D]=2\cdot D^2 \text{ and }0=[D,I]+[I,D]=2\cdot [D,I].$$
Having chosen a polarization
$P=\begin{pmatrix}
  D & 0 \\
  I & D \\
\end{pmatrix},$ the formula for
$\delta=\ad(P)$ is given by
$$
\delta
\begin{pmatrix}
  f & 0 \\
  i & f \\
\end{pmatrix}%
=
 \left[
\begin{pmatrix}
  D & 0 \\
  I & D \\
\end{pmatrix}%
,
\begin{pmatrix}
  f & 0 \\
  i & f \\
\end{pmatrix}%
\right]=
\begin{pmatrix}
  $[$D,f$]$ & 0 \\
  $[$D,i$]$+$[$f,I$]$ & $[$D,f$]$ \\
\end{pmatrix}.
$$

Now we look at the gauge equivalence.
First of all, the gauge group $G=\exp(\g^0)$
is the Lie group   consisting of
matrices of the form $e^A$, for any $A\in\g^0$. The gauge action of
$G$ on $\g$ is then determined by $e^{\ad(A)}\cdot B=\Ad(e^A)(B)=e^A B
 e^{-A}$. A computation shows that
$$ \exp{
\begin{pmatrix}
  f & 0 \\
  i & f \\
\end{pmatrix}%
}=
\begin{pmatrix}
  e^f & 0 \\
  x & e^f \\
\end{pmatrix},
\text{ where }
 x=\sum_{n\geq 1}\frac{1}{n!}\sum_{k+l=n-1} f^k\cdot i\cdot
f^l. $$ Then the gauge equivalence summarizes as
$$ e^A \begin{pmatrix}
  D & 0 \\
  I & D \\
\end{pmatrix}%
 e^{-A} =
\begin{pmatrix}
  e^f D e^{-f} & 0 \\
  e^f I e^{-f}+ [e^f D e^{-f},x e^{-f}]  & e^f D e^{-f} \\
\end{pmatrix}.$$
This concludes the example.
\end{example}

\begin{remark}\label{remark}
Let $N$ be a graded coalgebra over $S$.  
Then $\hom(N,N)$ will be a graded associative algebra 
and a Lie algebra with the bracket defined by the commutator.
The space $\coder(N)\subseteq \hom(N,N)$ will not, in general, form an associative algebra, but it is a Lie algebra.  
In particular, for any vector space $A$,  
$\coder(TA)$ is a graded Lie algebra.  
An $A_\infty$ structure on $A$ consists of an element
$D\in\coder(TA)$ satisfying $D^2=0$.  Thus one can say that an $A_\infty$ structure on $A$ is a choice of polarization $D\in\coder(TA)$.  Hence, if $(A,D)$ is an $A_\infty$ algebra,  $\coder(TA)$ carries a
differential $\delta:\coder(TA)\to \coder(TA)$ defined by
$$\delta(f):=[D,f]=D\circ f-(-1)^{|f|}f\circ D.$$
The complex ($\coder(TA),\delta$) is called the Hochschild cochain
complex of $A$.  Together with the bracket from $\hom(TA,TA)$, it is a differential graded Lie algebra that controls the deformations of the $A_\infty$ algebra $(A,D)$.  In order to make this statement precise, 
we recall the deformation theory of $A_\infty$ algebras (see for example \cite{FP}).  As a first observation, one may note that $\gamma$ is a solution to the Mauer-Cartan equation in the Hochschild dgLa if and only if $D+\gamma$ is another polarization in $\coder(TA)$; i.e., another $A_\infty$ structure, on $A$.  

\end{remark}

\subsection{Deformations of $A_\infty$ algebras}
Let $A$ be a graded vector space over a field $k$ and let $R$ be a
graded Artin local algebra with residue field $k$.  Let $m$ denote the maximal ideal
of $R$.  We have the decomposition $R \simeq R/m \oplus m \simeq k\oplus m$  and the projection $pr_k:R\to k$, hence the decomposition $A \otimes R \simeq A \oplus (A \otimes m)$ and the projection $pr_A:A \otimes R\to A$.  For definiteness, the reader may have the concrete example $R=k[t]/t^{l+1}$ in mind.  In this example, the maximal ideal is $m=tk[t]/t^{l+1}$,
$A \otimes R\simeq A+At+At^2+\cdots +A t^l$ (with the
tensor signs suppressed) and the natural projection $pr_A$ maps
$a_0+a_1t+a_2t^2+\cdots +a_lt^l\mapsto  a_0.$

Let $(A,D)$ be an $A_\infty$ algebra over $k$.  A \emph{deformation of $(A,D)$ over $R$} is
an $A_\infty$ algebra $(A\otimes R,D')$ over $R$ with the property that the projection
$$pr:T(A\otimes R)\simeq TA \otimes R \to TA $$ is
a morphism of $A_\infty$ algebras over $k$.  This means that $pr\circ D'=D \circ pr.$

Suppose that $D'$ is a deformation of $(A,D)$ over $R$.
Via any map $R\to S$, one can view $A\otimes R$ as an $S$ module
and $(A\otimes R,D')$ as a deformation of $(A,D)$ over $S$.

Let $\pi\in \hom(R\otimes R,R)$ denote the multiplication in $R$.
Let $D_R$ denote the $A_\infty$ structure $D \otimes \pi$ on
$A\otimes R.$  The $A_\infty$ algebra $(A\otimes R,D_R)$ is the
model for a trivial deformation of $(A,D)$.  That is, $(A\otimes R, D')$ 
is a \emph{trivial deformation}
if it is isomorphic to $(A\otimes R,D_R)$ as an
$A_\infty$ algebra.  This means that there is an
automorphism $$ \lambda:T(A\otimes R) \to T(A\otimes R)$$
satisfying $\lambda \circ D'= D_R \circ \lambda $.  Two deformations are equivalent if and only if they differ by a trivial one.

\subsection{Deformations of $A_\infty$ algebras with $\infty$ inner products}

\begin{definition}Let $A$ be a graded vector space over a field $k$.
  We define the graded Lie algebra $(\h=\oplus_i \h^i,[\,,\,])$ by
\begin{equation}
\h^i=\coder(TA)^{-i}\oplus \comap(T^{A}A,T^{A^{*}}A)^{1-i} \label{defpart1}\\
\end{equation}
and
\begin{multline} \label{defpart2}
  [(f,i),(g,j)]=([f,g], \delta_{f}(j)-(-1)^{|f||g|}\delta_{g}(i)) \\
  =(fg-(-1)^{|f||g|}gf, f^{A^*} j-(-1)^{|f||j|} j f^A -(-1)^{|f||g|} g^{A^*} i
    +(-1)^{|g|\cdot(|f|+|i|)} i g^A).
\end{multline}
\end{definition}
The skew-symmetry and Jacobi identity of $[\,,\,]$ are
straightforward to check after one notices that $\delta_f\circ
\delta_g -(-1)^{|f||g|} \delta_g\circ \delta_f = \delta_{ f\circ
g-(-1)^{|f||g|}g\circ f}$.

\begin{proposition}\label{prop3.1}
A pair $(D,I)\in \h$ is an $A_\infty$
  structure with $\infty$ inner product on $A$ if and only if
  $[(D,I),(D,I)]=0$.
\end{proposition}

\begin{proof}
  This is immediate:
$$ 0=[(D,I),(D,I)] \Leftrightarrow 0=[D,D]=2\cdot D^2 \text{ and }
  0=2\cdot \delta_D(I)=2(D^{A^*}\circ I - I \circ D^{A}).$$
The condition $D^2=0$ means that $D$ defines an $A_\infty$
structure on $A$ and the condition $D^{A^*}\circ I - I \circ
D^{A}=0$ means that $I$ defines a compatible $\infty$-inner
product.
\end{proof}

Now fix an $A_\infty$ structure together with an $\infty$ inner
product, which is to say, fix a pair $(D,I)\in \h$ with
$[(D,I),(D,I)]=0$.  Then, define $d:\h \to \h$ by
\begin{equation}\label{defpart3}
d(f,i)=[(D,I),(f,i)].
\end{equation}
The triple
$(\h,d,[\,,\,])$ is a differential graded Lie algebra.

\begin{definition}A \emph{deformation of an $A_\infty$ algebra
with $\infty$ inner product $(A,D,I)$ over $R$} is an
$A_\infty$ algebra over $R$ with $\infty$ inner product $(A\otimes R,D',
I')$, such that the projection
$$ pr:T(A\otimes R) \to TA $$
is a morphism of $A_\infty$ algebras over $k$ compatible with the
$\infty$-inner products. Compatibility with the $\infty$ inner product 
means that the following diagram of $A_\infty$-bimodule maps over
$k$ is commutative:
$$
\begin{diagram}
  \node{T^{A\otimes R} (A\otimes R)}\arrow{e,t}{\overline{pr}}
    \arrow{s,l}{pr_{k}\circ I'}
    \node{T^A (A\otimes R)}\arrow{s,r}{I^{pr}}\\
  \node{T^{(A\otimes R)^*} (A\otimes R)} \node{T^{A^*} (A\otimes R)}
    \arrow{w,b}{\widetilde{pr}}
\end{diagram}
$$
Here, the $\infty$-inner product $I'$ on $A\otimes R$ over $R$
induces an $\infty$-inner product on $A\otimes R$ over $k$ by
composing with the map induced by the projection $\hom_R(A\otimes
R,R)\to \hom_k(A\otimes R, k)$, $f\mapsto pr_{k}\circ f$.
\end{definition}

There is a natural extension of $I$ to an $\infty$-inner product
$I_R=I\otimes \pi$ on $(A\otimes R,D_R)$.

\begin{definition}\label{def_of_I}
We say that $(D',I')$ is a \emph{trivial deformation} of $(D,I)$
provided ``the triple $(A\otimes R,D',I')$ is isomorphic to
$(A\otimes R,D_R,I_R)$ as $A_\infty$ algebras with $\infty$ inner
products.''  That is, if there exist an
automorphism
\begin{gather*}
\lambda:  T(A\otimes R) \to T(A\otimes R)
\intertext{and a comap}
\rho: T^{A\otimes R}(A\otimes R) \to T^{(A\otimes
R)^*}(A\otimes R)
\end{gather*}
satisfying
\begin{enumerate}
\renewcommand{\theenumi}{\roman{enumi}}\renewcommand{\labelenumi}{(\theenumi)}
\item $\lambda \circ D'=D_R \circ \lambda$, \item $I'-
\widetilde{\lambda} \circ (I_R)^{\lambda} \circ
\overline{\lambda}= D'^{(A\otimes R)^*} \circ \rho + \rho \circ
D'^{A\otimes R}.$
\end{enumerate}
\end{definition}

It may be helpful to think of the second condition in
definition \ref{def_of_I} as saying
$I'$ equals $I_R$ under a change of coordinates (given by $\lambda$) up to a homotopy (given by $\rho$).  That is,
the following diagram commutes, up to a homotopy defined by $\rho \in \comap(T^{A\otimes R}(A\otimes R))$.

$$
\begin{diagram}
  \node{T^{A\otimes R} (A\otimes R)}\arrow{e,t}{\overline{\lambda}}
    \arrow{s,l}{I'} \node{T^{A\otimes R} (A\otimes R)}\arrow{s,r}{(I_{R})^\lambda}\\
  \node{T^{(A\otimes R)^*} (A\otimes R)}   \node{T^{(A\otimes R)^*} (A\otimes R)}
    \arrow{w,b}{\widetilde{\lambda}}
\end{diagram}
$$
Two deformations are equivalent if and only if they differ by a trivial one.

Now, the conclusion:
\begin{theorem}\label{mainthm}
Let $(A,D)$ be an $A_\infty$ algebra and let $I$ be an
$\infty$-inner product. Then the dg Lie algebra $(\h,d,[\,,\,])$
defined by equations \eqref{defpart1}, \eqref{defpart2} and
\eqref{defpart3} controls the deformations of the $A_\infty$
algebra with $\infty$ inner product $(A,D,I)$.
\end{theorem}

\begin{proof}
The content of this theorem is summarized in the following
two statements.  
\begin{itemize}
\item Deformations, over $R$, of the $(A,D,I)$ correspond to
  solutions to the Mauer-Cartan equation in $\h\otimes m$,
\item and equivalent deformations correspond to gauge
  equivalent solutions.
\end{itemize}

First we prove the first statement.
Let  $\alpha=(f,i)\in (\h\otimes m)^1$.
Observe that
\begin{align*}
  d\alpha + \frac{1}{2}[ \alpha,\alpha]
  &=[(D_R,I_R),(f,i)]+\frac{1}{2}\cdot[(f,i),(f,i)]\\
  &=\frac{1}{2}\cdot[(D_R+f,I_R+i),(D_R+f,I_R+i)].
\end{align*}
Then, proposition \ref{prop3.1} proves that $d\alpha + \frac{1}{2}[\alpha,\alpha]=0$
if and only if $(A\otimes R,D_R +f,I_R +i)$ is a
deformation of $(A,D,I)$.  It is immediate that any $(A\otimes R,D',I')$
that is a deformation of $(A,D,I)$ must satisfy $[(D',I'),(D',I')]=0 \in \h\otimes R$.  The fact that $pr:T(A\otimes R)\to TA$ is a map of $A_\infty$ algebras with $\infty$ inner products implies that $D'=D_R+f$ and $I'=I_R+i$ for some $(f,i)\in \h \otimes m$.

Now we prove the second statement.
Let
$\alpha=(f,i)\in (\h \otimes m)^0$.
The gauge action for $\h$ becomes
$$ e^{\ad(f,i)}\cdot (D_R,I_R)=\sum_{n\geq 0}
   \frac{\ad(f,i)^n}{n!}(D_R,I_R). $$
It follows from
$$ \delta_f\left(\delta_{\ad(f)^r(D_R)}((\delta_f)^s(i))\right)=
 \delta_{\ad(f)^{r+1}(D_R)}((\delta_f)^s(i))+
 \delta_{\ad(f)^r(D_R)}((\delta_f)^{s+1}(i)), $$
that $\ad(f,i)^n(D_R,I_R)$ is given by
$$  \left( \ad(f)^n(D_R) ,(\delta_f)^n(I_R)-
\sum_{k+l=n-1} \frac{n!}{k!(l+1)!}\cdot
\delta_{\ad(f)^k(D_R)}\circ(\delta_{f})^l(i)\right). $$ Now define
$\lambda^{-1}=e^{f}=\sum_{k\geq 0}\frac{1}{k!}f^k$ and $\rho=
\sum_{l\geq 0} \frac{-1}{(l+1)!} \cdot (\delta_{f})^l (i)$.  Then
for the automorphism $\lambda$ and the homotopy $\rho$, we have
\begin{align*}
\sum_{n\geq 0} &\frac{\ad(f,i)^n}{n!}(D_R,I_R)\\
&=
 \left( \sum_{n\geq 0} \frac{\ad(f)^n}{n!}(D_R),
    \sum_{n\geq 0} \frac{(\delta_f)^n}{n!}(I_R)+
     \delta_{\sum_{k\geq 0} \frac{\ad(f)^k}{k!}(D_R)}
       \left(\sum_{l\geq 0} \frac{-1}{(l+1)!} \cdot
         (\delta_{f})^l(i)\right)\right)\\
&=
 \left(\lambda^{-1} D_R \lambda,
   \widetilde{\lambda} (I_R)^{\lambda} \overline{\lambda}+
   \delta_{\lambda^{-1} D_R \lambda}(\rho)\right).
\end{align*}
This proves that $e^{\ad(f,i)}\cdot (D_R,I_R)$ is a trivial deformation of $(D,I)$.

It is not hard to see that every trivial deformation of $(D,I)$
arises from an element gauge equivalent to the identity. The
condition that the $A_\infty$ algebra map $\lambda:T(A\otimes
R)\to T(A\otimes R)$ is an automorphism implies that
$\lambda=e^{f}$ for some $f\in (\coder(TA)\otimes m)^0$. Also,
since $\rho=-i-\frac{1}{2}\delta_f (i)-\cdots$, the map $i\mapsto
\rho(i)=\sum_{l\geq 0} \frac{-1}{(l+1)!} \cdot (\delta_{f})^l (i)$
is invertible.  So one can obtain any homotopy $\rho$, by choosing
a suitable element $i= \sum_{m\geq 0}c_m\cdot
(\delta_f)^m(\rho)\in (\h\otimes m)^0$ with $\rho(i)=\rho$.
\end{proof}

\section{Moduli, infinitesimal deformations, and relationship to cyclic cohomology}

Let us return briefly to general deformation theory in order to
review the notions of infinitesimal deformations and moduli space.
Let $(\g,d,[\,,\,])$ be a differential graded Lie algebra and
assume that $\ker(d)/\im(d)=:H(\g)=\oplus_{i=-m}^m H^i(\g)$ is
finite dimensional.  Consider the (graded version of the)
ring of dual numbers $R=k[t_{-m}, \ldots, t_m]/t_it_j$.  Here
$\deg(t_i)=i-1$ and the maximal ideal of $R$ is $m=\oplus_i t_i
R.$

From a solution
 $\sum(\gamma_j \otimes t_j)\in (\g \otimes m)^1$ to the Mauer-Cartan equation,
one may produce the map $d +\sum t_j\ad(\gamma_j):\g \otimes k[t_{-m},\ldots, t_m] \to \g \otimes k[t_{-m},\ldots, t_m]$
which satisfies
$$\left(d +\sum t_j\ad(\gamma_j)\right)^2=0 \text{ modulo }t_it_j.$$
One refers to $\gamma=\sum \gamma_j$ as an infinitesimal deformation.
One can readily check that
$$Def_\g(R)=\ker(d)/\im(d)=H(\g).$$
Schlessinger's theorem \cite{S} implies that $\Def_\g$ is
prorepresentable.  That is, there exists a projective limit of
(graded) local Artin rings $\O$ and an equivalence of the functors
$$\Def_\g(\cdot)\simeq \hom(\O,\cdot).$$
In the case that $\O=\O_\M$ is the ring of local functions at the
base point of a pointed $\mathbb{Z}$ graded space $\M$, then $\M$
is the local moduli space for $\Def_\g$.
Denote the base point of $\M$ by $0$.  One can check that $$T_0(\M)\simeq\hom(\O_\M,R).$$  It follows that the
graded tangent space to the moduli space at the base point 
is isomorphic to
the cohomology of $(\g,d)$:
$$T_0(\M)\simeq H(\g).$$

Now, let $(A,D)$ be an $A_\infty$ algebra and let $I$ be an $\infty$ inner product on $(A,D)$.
Theorem \ref{mainthm} says that the dg Lie algebra controlling deformations
of  $(A,D,I)$ is
\begin{gather*}
  \h=\coder(TA)\oplus \comap(T^{A}A,T^{A^{*}}A)
\intertext{with bracket}
  [(f,i),(g,j)] =([f,g], \delta_{f}(j)-(-1)^{|f||g|}\delta_{g}(i))
\intertext{and a differential} d(f,i) = [(D,I),(f,i)].
\end{gather*}
Thus follows the expected infinitesimal statement:
\begin{corollary}
The tangent space to the moduli space of $A_\infty$ structures
with $\infty$  inner products is isomorphic to $H(\h).$
\end{corollary}

As a final remark, we mention some connections between the
cohomology $H(\h)$ and a couple of its cousins.  If $(A,D,I)$ is
an $A_\infty$ algebra with $\infty$-inner product, we have the
Hochschild dgLa $(\coder(TA),\delta,[\,,\,])$ and the sub dgLa of
cyclic Hochschild cochains $\coder(TA)_{\mathrm{Cyclic}}$, defined
by
$$\coder(TA)_{\mathrm{Cyclic}}=\{f \in \coder(TA): \delta_f(I)=
0\}.$$ Note from appendix A that if $I$ consists of an ordinary
symmetric inner product $I=\<\,,\,\>$, then the condition
$\delta_f(I)=f^{A^*}\circ I-I \circ f^A=0$ is equivalent to
$$ \<f(a_1,...,a_n),a_{n+1}\>=\pm \<a_1,f(a_2,...,a_{n+1})\>. $$

We have the following maps of differential graded Lie algebras.
\begin{equation}\label{maps}
(\coder(TA)_{\mathrm{Cyclic}},\delta,[\,,\,]) \longrightarrow
(\h,d,[\,,\,]) \text{ and } (\h,d,[\,,\,])
\longrightarrow (\coder(TA),\delta,[\,,\,]).
\end{equation}
 The first map is the injection $f\mapsto (f,0)\in \h,$ which is
 is a cochain map
\begin{equation*}\label{cyclic1}
d(f,0)=([D,f],\pm(f^{A^*}\circ I-I \circ f^{A}))= (\delta f,0),
\end{equation*}
because elements of the domain are cyclic.  The induced map in
cohomology describes a statement from \cite{P}, namely that the
first order deformations of $D$ compatible with the inner product
are classified by cyclic cohomology.  We do not know under what conditions
the map $f \mapsto (f,0)\in \h$ induces an isomorphism in cohomology.
 The second map in (\ref{maps}) is simply the projection $\coder(TA)\oplus
\comap(T^{A}A,T^{A^{*}}A) \to \coder(TA)$ and the induced map in
cohomology describes the simple statement that any infinitesimal
deformation of the pair $(D,I)$ gives an infinitesimal deformation
of $D$.

\appendix
\section{Explicit formulas of $\delta_f(i)$}

Let $f\in \coder(TA)$ and $i\in \comap(T^{A}A,T^{A^{*}}A)$. We
want to describe the term $\delta_f(i)=f^{A^*}\circ
i-(-1)^{|f||i|} \cdot i\circ f^{A} \in \comap(T^{A}A,T^{A^{*}}A)$
more explicitly. Here, $f:\oplus_{k\geq 1} A^{\otimes k}\to A$ and
$i:\oplus_{k,l\geq 0} A^{\otimes k}\otimes A \otimes A^{\otimes l}
\otimes A \to S$ have the components
$$
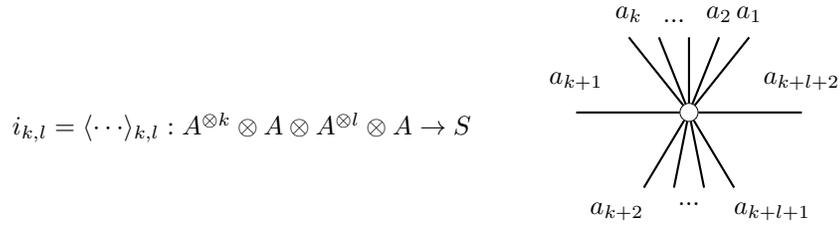
\begin{pspicture}(-7,0)(6,4)
  \rput[tl](-4,2){$f_k:A^{\otimes k}\to A$}
 \psline(2,2)(1.2,3)
 \psline(2,2)(1.6,3)
 \psline(2,2)(2,3)
 \psline(2,2)(2.4,3)
 \psline(2,2)(2.8,3)
 \psline(2,2)(2,1)
 \psdots[dotstyle=*,dotscale=2](2,2)
 \rput[b](2.8,3.2){$a_1$}  \rput[b](2.4,3.2){$a_2$}
 \rput[b](1.8,3.2){$...$} \rput[b](1.2,3.2){$a_k$}
 \rput[t](2,.8){$f_k(a_1,...,a_k)$}
\end{pspicture}
$$
$$
\begin{pspicture}(-7,0)(6,4)
  \rput[tl](-7,2){$i_{k,l}=\<\cdots \>_{k,l}:A^{\otimes k}\otimes A \otimes
A^{\otimes l} \otimes A \to S$}
 \psline(.5,2)(3.5,2)
 \psline(2,2)(1.2,3)
 \psline(2,2)(1.6,3)
 \psline(2,2)(2,3)
 \psline(2,2)(2.4,3)
 \psline(2,2)(2.8,3)
 \psline(2,2)(1.4,1)
 \psline(2,2)(1.8,1)
 \psline(2,2)(2.2,1)
 \psline(2,2)(2.6,1)
 \psdots[dotstyle=o,dotscale=2](2,2)
 \rput[b](2.8,3.2){$a_1$}  \rput[b](2.4,3.2){$a_2$}
 \rput[b](1.8,3.2){$...$} \rput[b](1.2,3.2){$a_k$}
 \rput[b](.5,2.3){$a_{k+1}$}
 \rput[b](3.5,2.3){$a_{k+l+2}$}
 \rput[tr](1.4,.8){$a_{k+2}$} \rput[t](2,.8){$...$}
 \rput[tl](2.6,.8){$a_{k+l+1}$}
\end{pspicture}
$$
By convention, the inputs are always inserted using the
counterclockwise direction.
Then $f^{A^*}\circ i-(-1)^{|f||i|}\cdot i\circ f^{A}$ is given by
inserting $f$ into $i$ in all possible combinations.
$$ \pm \,\,\,
     \pstree[treemode=R, levelsep=1cm, treesep=0.3cm]{\Tp}
    { \pstree[levelsep=0cm]{
        \Tr{\Tc{4pt} }}
  { \pstree[treemode=U, levelsep=0.8cm]{\Tn}
     {\Tp \Tp \Tp \pstree{\Tc*{4pt}}{\Tp \Tp \Tp \Tp} \Tp}
    \pstree[treemode=R, levelsep=1cm]{\Tn}
     {\Tp}
    \pstree[treemode=D, levelsep=0.8cm]{\Tn}
     {\Tp \Tp \Tp \Tp \Tp}
  }}
 \,\,\, \pm \,\,\,
     \pstree[treemode=L, levelsep=1cm, treesep=0.3cm]{\Tp}
    { \pstree[levelsep=0cm]{
        \Tr{\Tc{4pt} }}
  { \pstree[treemode=U, levelsep=0.8cm]{\Tn}
     {\Tp \Tp \Tp \Tp \Tp}
    \pstree[treemode=L, levelsep=0.7cm]{\Tn}
     {\pstree{\Tc*{4pt}}{\Tp \Tp \Tp \Tp \Tp}}
    \pstree[treemode=D, levelsep=0.8cm]{\Tn}
     {\Tp \Tp \Tp \Tp \Tp}
  }}
 \,\,\, \pm \,\,\,
     \pstree[treemode=R, levelsep=1cm, treesep=0.3cm]{\Tp}
    { \pstree[levelsep=0cm]{
        \Tr{\Tc{4pt} }}
  { \pstree[treemode=U, levelsep=0.8cm]{\Tn}
     {\Tp \Tp \Tp \Tp \Tp}
    \pstree[treemode=R, levelsep=1cm]{\Tn}
     {\Tp}
    \pstree[treemode=D, levelsep=0.8cm]{\Tn}
     {\Tp \Tp \Tp \pstree{\Tc*{4pt}}{\Tp \Tp \Tp \Tp} \Tp}
  }}
 \,\,\, \pm \,\,\,
     \pstree[treemode=R, levelsep=1cm, treesep=0.3cm]{\Tp}
    { \pstree[levelsep=0cm]{
        \Tr{\Tc{4pt} }}
  { \pstree[treemode=U, levelsep=0.8cm]{\Tn}
     {\Tp \Tp \Tp \Tp \Tp}
    \pstree[treemode=R, levelsep=0.7cm]{\Tn}
     {\pstree{\Tc*{4pt}}{\Tp \Tp \Tp \Tp \Tp}}
    \pstree[treemode=D, levelsep=0.8cm]{\Tn}
     {\Tp \Tp \Tp \Tp \Tp}
  }}
$$

First, here are some examples of how these diagrams are to be
read.

$\<a,b,c,d\>_{2,0}$\\
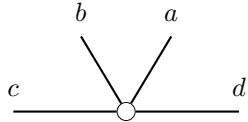
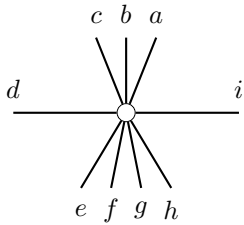
$$
\begin{pspicture}(0,0)(4,3)
 \psline(.5,1)(3.5,1)
 \psline(2,1)(1.4,2)
 \psline(2,1)(2.6,2)
 \psdots[dotstyle=o,dotscale=2](2,1)
 \rput[b](2.6,2.2){$a$}  \rput[b](1.4,2.2){$b$}
 \rput[b](0.5,1.2){$c$}  \rput[b](3.5,1.2){$d$}
\end{pspicture}
$$
$\<a,b,c,d,e,f,g,h,i\>_{3,4}$\\
$$
\begin{pspicture}(0,0)(4,4)
 \psline(.5,2)(3.5,2)
 \psline(2,2)(1.6,3)
 \psline(2,2)(2,3)
 \psline(2,2)(2.4,3)
 \psline(2,2)(1.4,1)
 \psline(2,2)(1.8,1)
 \psline(2,2)(2.2,1)
 \psline(2,2)(2.6,1)
 \psdots[dotstyle=o,dotscale=2](2,2)
 \rput[b](2,3.2){$b$}  \rput[b](2.4,3.2){$a$}
 \rput[b](1.6,3.2){$c$}
 \rput[b](.5,2.2){$d$} \rput[b](3.5,2.2){$i$}
 \rput(1.4,.7){$e$} \rput(1.8,.7){$f$}
 \rput(2.2,.7){$g$} \rput(2.6,.7){$h$}
\end{pspicture}
$$
$\<f_{2}(f_{2}(b,c),f_{2}(d,e)),f_{2}(f,a)\>_{0,0}$\\
$$
\begin{pspicture}(0,0)(4,4)
 \psline(.5,2)(3.5,2)
 \psline(3,2)(3.5,3)
 \psline(1,2)(.5,3)
 \psline(1.66,2)(.74,1)
 \psline(1.2,1.5)(1.2,.8)
 \psdots[dotstyle=*,dotscale=2](1.66,2)
 \psdots[dotstyle=*,dotscale=2](1,2)
 \psdots[dotstyle=*,dotscale=2](1.2,1.5)
 \psdots[dotstyle=*,dotscale=2](3,2)
 \psdots[dotstyle=o,dotscale=2](2.33,2)
 \rput[b](.2,2){$c$}    \rput[b](3.8,2){$f$}
 \rput[b](3.5,3.2){$a$} \rput[b](0.5,3.2){$b$}
 \rput[b](.5,.8){$d$}   \rput[b](1.2,.5){$e$}
\end{pspicture}
$$
$\<a,b,f_{3}(c,d,f_{2}(e,f)),g,f_{2}(h,i))\>_{1,2}$\\
$$
\begin{pspicture}(0,-.8)(4,4)
 \psline(.5,2)(3.5,2)
 \psline(2,2)(2,3)
 \psline(2,2)(.8,.8)
 \psline(2,2)(2.6,1)
 \psline(2.8,2)(3.4,1.4)
 \psline(1.3,1.3)(.6,1.2)
 \psline(1.3,1.3)(1.4,.6)
 \psline(1.4,.6)(1.2,.2)
 \psline(1.4,.6)(1.7,.2)
 \psdots[dotstyle=o,dotscale=2](2,2)
 \psdots[dotstyle=*,dotscale=2](2.8,2)
 \psdots[dotstyle=*,dotscale=2](1.3,1.3)
 \psdots[dotstyle=*,dotscale=2](1.4,.6)
 \rput[b](2,3.2){$a$} \rput[b](.5,2.2){$b$}
 \rput(2.6,.7){$g$}   \rput[b](3.5,2.2){$i$}
 \rput[b](3.6,1.2){$h$}
 \rput[b](.4,1.1){$c$} \rput[b](.6,.4){$d$}
 \rput[b](1.1,-.2){$e$} \rput[b](1.7,-.2){$f$}
\end{pspicture}
$$
$\<c,f_{2}(d,e),f_{2}(f_{2}(f,g),h),i,f_{4}(j,k,a,b)\>_{2,1}$\\
$$
\begin{pspicture}(0,0)(4,4)
 \psline(.5,2)(3.5,2)
 \psline(3,2)(3.5,3)
 \psline(1,2)(.5,3)
 \psline(1.66,2)(.74,1)
 \psline(3,2)(3.6,2.5)
 \psline(2.33,2)(2,1)
 \psline(3,2)(3.6,1.2)
 \psline(2.33,2)(2.6,3)
 \psline(2.33,2)(2,2.5)
 \psline(2,2.5)(2,3)
 \psline(2,2.5)(1.5,3)
 \psdots[dotstyle=*,dotscale=2](1.66,2)
 \psdots[dotstyle=*,dotscale=2](1,2)
 \psdots[dotstyle=*,dotscale=2](3,2)
 \psdots[dotstyle=*,dotscale=2](2,2.5)
 \psdots[dotstyle=o,dotscale=2](2.33,2)
 \rput[b](.2,1.9){$g$}  \rput[b](3.8,1.9){$k$}
 \rput[b](3.5,3.2){$b$} \rput[b](0.5,3.2){$f$}
 \rput[b](.5,.8){$h$}   \rput[b](3.8,2.5){$a$}
 \rput[b](2,.6){$i$}    \rput[b](3.8,1){$j$}
 \rput[b](2.6,3.2){$c$} \rput[b](2,3.2){$d$}
 \rput[b](1.5,3.2){$e$}
\end{pspicture}
$$

Here are the terms of $\delta_f(i)=f^{A^*}\circ
i-(-1)^{|f||i|}\cdot i\circ f^{A}$ up to sign, when they are being
applied to elements from
$A^{\otimes k}\otimes A \otimes A^{\otimes l}\otimes A$:\\
\underline{ $k=0$, $l=0$:}\quad
$$\<f_1(a),b\>_{0,0}\pm\<a,f_1(b)\>_{0,0} $$
$$
\begin{pspicture}(4,0)(9,2)

 \psline(4.3,1)(5.7,1)
 \psdots[dotstyle=*,dotscale=1.4](4.8,1)
 \psdots[dotstyle=o,dotscale=1.4](5.3,1)
 \rput[b](4.3,1.1){$a$} \rput[b](5.7,1.1){$b$}

 \rput(6.5,1.1){$\pm$}

 \psline(7.3,1)(8.8,1)
 \psdots[dotstyle=*,dotscale=1.4](8.3,1)
 \psdots[dotstyle=o,dotscale=1.4](7.8,1)
 \rput[b](7.3,1.1){$a$} \rput[b](8.8,1.1){$b$}
\end{pspicture}
$$
\underline{$k=1$, $l=0$:}
\begin{eqnarray*}
\<f_1(a),b,c\>_{1,0}\pm
\<a,f_1(b),c\>_{1,0}\pm \<a,b,f_1(c)\>_{1,0}\pm & & \\
\<f_2(a, b),c\>_{0,0}\pm \<b,f_2(c,a)\>_{0,0} &&
\end{eqnarray*}
$$
\begin{pspicture}(0,0)(10,2.5)
 \psline(1.3,1)(2.7,1)
 \psline(2,1)(2,1.8)
 \psdots[dotstyle=o,dotscale=1.4](2,1)
 \psdots[dotstyle=*,dotscale=1.4](2,1.4)
 \rput[b](1.3,1.1){$b$} \rput[b](2.7,1.1){$c$} \rput[b](2.2,1.6){$a$}
 \rput(3.4,1.1){$\pm$}

 \psline(4.3,1)(5.7,1)
 \psline(5,1)(5,1.8)
 \psdots[dotstyle=o,dotscale=1.4](5,1)
 \psdots[dotstyle=*,dotscale=1.4](4.6,1)
 \rput[b](4.3,1.1){$b$} \rput[b](5.7,1.1){$c$} \rput[b](5.2,1.6){$a$}
 \rput(6.4,1.1){$\pm$}

 \psline(7.3,1)(8.7,1)
 \psline(8,1)(8,1.8)
 \psdots[dotstyle=o,dotscale=1.4](8,1)
 \psdots[dotstyle=*,dotscale=1.4](8.4,1)
 \rput[b](7.3,1.1){$b$} \rput[b](8.7,1.1){$c$} \rput[b](8.2,1.6){$a$}
 \rput(9.4,1.1){$\pm$}
\end{pspicture}
$$
$$
\begin{pspicture}(0,0)(10,2.5)
 \psline(4.3,1)(5.7,1)
 \psline(4.8,1)(4.6,1.8)
 \psdots[dotstyle=*,dotscale=1.4](4.8,1)
 \psdots[dotstyle=o,dotscale=1.4](5.3,1)
 \rput[b](4.3,1.1){$b$} \rput[b](5.7,1.1){$c$} \rput[b](4.8,1.8){$a$}

 \rput(6.5,1.1){$\pm$}

 \psline(7.3,1)(8.8,1)
 \psline(8.3,1)(8.5,1.8)
 \psdots[dotstyle=*,dotscale=1.4](8.3,1)
 \psdots[dotstyle=o,dotscale=1.4](7.8,1)
 \rput[b](7.3,1.1){$b$} \rput[b](8.8,1.1){$c$} \rput[b](8.3,1.8){$a$}
\end{pspicture}
$$
\underline{$k=0$, $l=1$:}
\begin{eqnarray*}
\<f_1(a),b,c\>_{0,1}\pm
\<a,f_1(b),c\>_{0,1}\pm \<a,b,f_1(c)\>_{0,1}\pm & & \\
\<f_2(a, b),c\>_{0,0}\pm \<a,f_2(b,c)\>_{0,0} & &
\end{eqnarray*}
$$
\begin{pspicture}(0,-.2)(10,1.8)
 \psline(1.3,1)(2.7,1)
 \psline(2,1)(2,0.2)
 \psdots[dotstyle=o,dotscale=1.4](2,1)
 \psdots[dotstyle=*,dotscale=1.4](1.6,1)
 \rput[b](1.3,1.1){$a$} \rput[b](2.7,1.1){$c$} \rput[b](2.2,0.2){$b$}
 \rput(3.4,1.1){$\pm$}

 \psline(4.3,1)(5.7,1)
 \psline(5,1)(5,0.2)
 \psdots[dotstyle=o,dotscale=1.4](5,1)
 \psdots[dotstyle=*,dotscale=1.4](5,0.6)
 \rput[b](4.3,1.1){$a$} \rput[b](5.7,1.1){$c$} \rput[b](5.2,0.2){$b$}
 \rput(6.4,1.1){$\pm$}

 \psline(7.3,1)(8.7,1)
 \psline(8,1)(8,0.2)
 \psdots[dotstyle=o,dotscale=1.4](8,1)
 \psdots[dotstyle=*,dotscale=1.4](8.4,1)
 \rput[b](7.3,1.1){$a$} \rput[b](8.7,1.1){$c$} \rput[b](8.2,0.2){$b$}
 \rput(9.4,1.1){$\pm$}
\end{pspicture}
$$
$$
\begin{pspicture}(0,-.2)(10,1.8)

 \psline(4.3,1)(5.7,1)
 \psline(4.8,1)(4.6,0.2)
 \psdots[dotstyle=*,dotscale=1.4](4.8,1)
 \psdots[dotstyle=o,dotscale=1.4](5.3,1)
 \rput[b](4.3,1.1){$a$} \rput[b](5.7,1.1){$c$} \rput[b](4.8,0.2){$b$}

 \rput(6.5,1.1){$\pm$}

 \psline(7.3,1)(8.8,1)
 \psline(8.3,1)(8.5,0.2)
 \psdots[dotstyle=*,dotscale=1.4](8.3,1)
 \psdots[dotstyle=o,dotscale=1.4](7.8,1)
 \rput[b](7.3,1.1){$a$} \rput[b](8.8,1.1){$c$} \rput[b](8.3,0.2){$b$}
\end{pspicture}
$$
\underline{$k=2$, $l=0$:}
\begin{eqnarray*}
\<f_1(a),b,c,d\>_{2,0}\pm \<a,f_1(b),c,d\>_{2,0}\pm && \\
\<a,b,f_1(c),d\>_{2,0}\pm \<a,b,c,f_1(d)\>_{2,0}\pm && \\
\<f_2(a,b),c,d\>_{1,0}\pm \<a,f_2(b,c),d\>_{1,0}\pm \<b,c,f_2(d,a)\>_{1,0}\pm &&\\
\<f_3(a,b,c),d\>_{0,0}\pm \<c,f_3(d,a,b)\>_{0,0}& &
\end{eqnarray*}
Note that for example the term $\<a,b,f_2(c,d)\>_{2,0}$ does not
appear, because $c$ and $d$ are the two special elements of
$a\otimes b\otimes c\otimes d\in A^{\otimes 2}\otimes A\otimes
A^{\otimes 0}\otimes A$, which are put on the horizontal line of
the diagram. The two special elements from $A^{\otimes k}\otimes
A\otimes A^{\otimes l}\otimes A$ can never be inside any $f_n$.
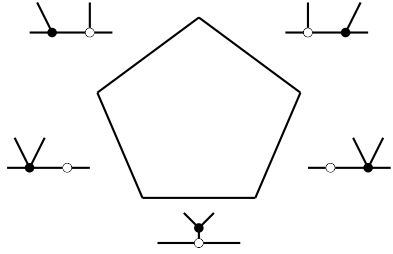
$$
\begin{pspicture}(3,0)(10,6)

 \psline(5.2,3.2)(6.55,4.2) \psline(6.55,4.2)(7.9,3.2)
 \psline(7.9,3.2)(7.3,1.8) \psline(7.3,1.8)(5.8,1.8)  \psline(5.8,1.8)(5.2,3.2)

 \psline(6,1.2)(7.1,1.2)
 \psline(6.55,1.4)(6.75,1.6) \psline(6.55,1.4)(6.35,1.6)
 \psdots[dotstyle=*,dotscale=1](6.55,1.4)
 \psline(6.55,1.2)(6.55,1.4) \psdots[dotstyle=o,dotscale=1](6.55,1.2)

 \psline(4,2.2)(5.1,2.2)
 \psline(4.3,2.2)(4.5,2.6) \psdots[dotstyle=*,dotscale=1](4.3,2.2)
 \psline(4.3,2.2)(4.1,2.6) \psdots[dotstyle=o,dotscale=1](4.8,2.2)

 \psline(8,2.2)(9.1,2.2)
 \psline(8.8,2.2)(8.6,2.6) \psdots[dotstyle=o,dotscale=1](8.3,2.2)
 \psline(8.8,2.2)(9,2.6) \psdots[dotstyle=*,dotscale=1](8.8,2.2)

 \psline(4.3,4)(5.4,4)
 \psline(4.6,4)(4.4,4.4) \psdots[dotstyle=*,dotscale=1](4.6,4)
 \psline(5.1,4)(5.1,4.4) \psdots[dotstyle=o,dotscale=1](5.1,4)

 \psline(7.7,4)(8.8,4)
 \psline(8,4)(8,4.4) \psdots[dotstyle=o,dotscale=1](8,4)
 \psline(8.5,4)(8.7,4.4) \psdots[dotstyle=*,dotscale=1](8.5,4)
\end{pspicture}
$$
\underline{$k=0$, $l=2$:}
\begin{eqnarray*}
\<f_1(a),b,c,d\>_{0,2}\pm \<a,f_1(b),c,d\>_{0,2}\pm && \\
\<a,b,f_1(c),d\>_{0,2}\pm \<a,b,c,f_1(d)\>_{0,2}\pm && \\
\<f_2(a,b),c,d\>_{0,1}\pm \<a,f_2(b,c),d\>_{0,1}\pm \<a,b,f_2(c,d)\>_{0,1}\pm &&\\
\<f_3(a,b,c),d\>_{0,0}\pm \<a,f_3(b,c,d)\>_{0,0} & &
\end{eqnarray*}
The special elements are $a$ and $d$ from $a\otimes b\otimes
c\otimes d\in A^{\otimes 0}\otimes A\otimes A^{\otimes 2}\otimes
A$.
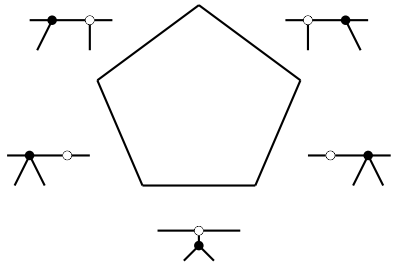
$$
\begin{pspicture}(3,0)(10,6)

 \psline(5.2,3.2)(6.55,4.2) \psline(6.55,4.2)(7.9,3.2)
 \psline(7.9,3.2)(7.3,1.8) \psline(7.3,1.8)(5.8,1.8)  \psline(5.8,1.8)(5.2,3.2)

 \psline(6,1.2)(7.1,1.2)
 \psline(6.55,1)(6.75,.8) \psline(6.55,1)(6.35,.8)
 \psdots[dotstyle=*,dotscale=1](6.55,1)
 \psline(6.55,1.2)(6.55,1) \psdots[dotstyle=o,dotscale=1](6.55,1.2)

 \psline(4,2.2)(5.1,2.2)
 \psline(4.3,2.2)(4.5,1.8) \psdots[dotstyle=*,dotscale=1](4.3,2.2)
 \psline(4.3,2.2)(4.1,1.8) \psdots[dotstyle=o,dotscale=1](4.8,2.2)

 \psline(8,2.2)(9.1,2.2)
 \psline(8.8,2.2)(8.6,1.8) \psdots[dotstyle=o,dotscale=1](8.3,2.2)
 \psline(8.8,2.2)(9,1.8) \psdots[dotstyle=*,dotscale=1](8.8,2.2)

 \psline(4.3,4)(5.4,4)
 \psline(4.6,4)(4.4,3.6) \psdots[dotstyle=*,dotscale=1](4.6,4)
 \psline(5.1,4)(5.1,3.6) \psdots[dotstyle=o,dotscale=1](5.1,4)

 \psline(7.7,4)(8.8,4)
 \psline(8,4)(8,3.6) \psdots[dotstyle=o,dotscale=1](8,4)
 \psline(8.5,4)(8.7,3.6) \psdots[dotstyle=*,dotscale=1](8.5,4)
\end{pspicture}
$$
\underline{$k=1$, $l=1$:}
\begin{eqnarray*}
\<f_1(a),b,c,d\>_{1,1}\pm \<a,f_1(b),c,d\>_{1,1}\pm && \\
\<a,b,f_1(c),d\>_{1,1}\pm \<a,b,c,f_1(d)\>_{1,1}\pm && \\
\<f_2(a,b),c,d\>_{0,1}\pm \<b,c,f_2(d,a)\>_{0,1}\pm &&\\
\<a,f_2(b,c),d\>_{1,0}\pm \<a,b,f_2(c,d)\>_{1,0}\pm &&\\
\<f_3(a,b,c),d\>_{0,0}\pm \<b,f_3(c,d,a)\>_{0,0}& &
\end{eqnarray*}
The special elements are $b$ and $d$ from $a\otimes b\otimes
c\otimes d\in A^{\otimes 1}\otimes A\otimes A^{\otimes 1} \otimes
A$.
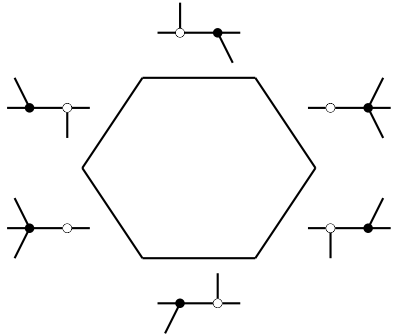
$$
\begin{pspicture}(3,0)(10,6)

 \psline(5,3)(5.8,4.2) \psline(5.8,4.2)(7.3,4.2)  \psline(7.3,4.2)(8.1,3)
 \psline(8.1,3)(7.3,1.8) \psline(7.3,1.8)(5.8,1.8)  \psline(5.8,1.8)(5,3)

 \psline(6,1.2)(7.1,1.2)
 \psline(6.3,1.2)(6.1,.8)  \psdots[dotstyle=*,dotscale=1](6.3,1.2)
 \psline(6.8,1.2)(6.8,1.6) \psdots[dotstyle=o,dotscale=1](6.8,1.2)

 \psline(6,4.8)(7.1,4.8)
 \psline(6.3,4.8)(6.3,5.2) \psdots[dotstyle=o,dotscale=1](6.3,4.8)
 \psline(6.8,4.8)(7,4.4) \psdots[dotstyle=*,dotscale=1](6.8,4.8)

 \psline(4,2.2)(5.1,2.2)
 \psline(4.3,2.2)(4.1,2.6) \psdots[dotstyle=*,dotscale=1](4.3,2.2)
 \psline(4.3,2.2)(4.1,1.8) \psdots[dotstyle=o,dotscale=1](4.8,2.2)

 \psline(8,2.2)(9.1,2.2)
 \psline(8.3,2.2)(8.3,1.8) \psdots[dotstyle=o,dotscale=1](8.3,2.2)
 \psline(8.8,2.2)(9,2.6) \psdots[dotstyle=*,dotscale=1](8.8,2.2)

 \psline(4,3.8)(5.1,3.8)
 \psline(4.3,3.8)(4.1,4.2) \psdots[dotstyle=*,dotscale=1](4.3,3.8)
 \psline(4.8,3.8)(4.8,3.4) \psdots[dotstyle=o,dotscale=1](4.8,3.8)

 \psline(8,3.8)(9.1,3.8)
 \psline(8.8,3.8)(9,4.2) \psdots[dotstyle=o,dotscale=1](8.3,3.8)
 \psline(8.8,3.8)(9,3.4) \psdots[dotstyle=*,dotscale=1](8.8,3.8)
\end{pspicture}
$$
\underline{$i=\<\, ,\, \>_{0,0}$ for any $k$, $l$:} Assume that
$i=\<\, ,\, \>_{0,0}$ has only lowest component, but $f$ has all
higher components. We apply $f^{A^*}\circ i-(-1)^{|f||i|}\cdot
i\circ f^{A}$ to the element
$$ a_1\otimes...\otimes a_k \otimes a_{k+1} \otimes a_{k+2} \otimes
 ... \otimes a_{k+l+1} \otimes a_{k+l+2} \in A^{\otimes k}\otimes
 A\otimes A^{\otimes l}\otimes A $$
to get
$$ \<f(a_1,...,a_{k+l+1}),a_{k+l+2}\>_{0,0}\pm
   \<a_{k+1},f(a_{k+2},...,a_{k+l+2},a_1,...,a_k)\>_{0,0} $$
$$
\begin{pspicture}(4,0)(10.5,2)

 \psline(4.3,1)(5.7,1)
 \psline(4.8,1)(4.3,1.8)  \psline(4.8,1)(4.3,1.4)
 \psline(4.8,1)(4.3,0.2)  \psline(4.8,1)(4.3,0.6)
 \psdots[dotstyle=*,dotscale=1.4](4.8,1)
 \psdots[dotstyle=o,dotscale=1.4](5.3,1)

 \rput(6.5,1.1){$\pm$}

 \psline(7.3,1)(8.8,1)
 \psline(8.3,1)(8.8,1.8)  \psline(8.3,1)(8.8,1.4)
 \psline(8.3,1)(8.8,0.2)  \psline(8.3,1)(8.8,0.6)
 \psdots[dotstyle=*,dotscale=1.4](8.3,1)
 \psdots[dotstyle=o,dotscale=1.4](7.8,1)
\end{pspicture}
$$

%
%

\bigskip

\end{document}